\documentclass{amsart}
\usepackage{amssymb, amsmath, color}
\usepackage[latin1]{inputenc}
\usepackage{color}

\newtheorem{theorem}{Theorem}
\newtheorem{lemma}[theorem]{Lemma}
\newtheorem{remark}[theorem]{Remark}
\newtheorem{definition}[theorem]{Definition}

\newtheorem{example}[theorem]{Example}
\newcommand{\grad}{\mbox{\rm grad}\,}

\begin{document}

\title[Locally homogeneous Lorentzian gradient Ricci solitons]
{The structure of the Ricci tensor on locally homogeneous Lorentzian 
gradient Ricci solitons}
\author{M. Brozos-V\'{a}zquez \, E. Garc\'{i}a-R\'{i}o  \, S. Gavino-Fern\'{a}ndez \, P. Gilkey}
\address{MBV: Departmento de Matemáticas, Escola Polit\'ecnica Superior, Universidade da Coru\~na, Spain}
\email{miguel.brozos.vazquez@udc.gal}
\address{PG: Mathematics Department, \; University of Oregon, \;\;
  Eugene \; OR 97403, \; USA}
\email{gilkey@uoregon.edu}
\address{EGR-SGF: Faculty of Mathematics,
University of Santiago de Compostela,
15782 Santiago de Compostela, Spain}
\email{eduardo.garcia.rio@usc.es $\,\,$ sandra.gavino@usc.es}
\thanks{Partially supported by projects EM2014/009, GRC2013-045 and MTM2013-41335-P with FEDER funds (Spain).}

\begin{abstract}
We describe the structure of the Ricci tensor on a locally
homogeneous Lorentzian gradient Ricci soliton. In 
the non-steady case, we show the soliton is rigid in dimensions three and four.
In the steady case, we give a complete classification in dimension three.
\end{abstract}

\maketitle

\section{Introduction}
Let $(M,g)$ be a Lorentzian manifold of dimension $n+2$ for $n\ge1$, let
$\rho$ be the {\it Ricci tensor}. Let $\operatorname{Ric}$
be the {\it Ricci operator}; $\rho(X,Y)=g(\operatorname{Ric}X,Y)$. 
If $f\in C^\infty(M)$, let
$\operatorname{Hess}_f$ be the Hessian; $f$ is often called the
{\it potential function}. Then 
$$\operatorname{Hess}_f(X,Y)=(\nabla_X df)(Y)=XY(f)-(\nabla_XY)(f).$$ 
Let $\nabla f$ be the vector field 
dual to the exterior derivative $df$ of $f$; this will also be denoted by $\grad \{ f\}$ for notational clarity when convenient. 
The {\it Hessian operator} 
$$\mathcal{H}_f(X):=\nabla_X(\nabla f)\text{ satisfies }
\operatorname{Hess}(X,Y)=g(\mathcal{H}_fX,Y)\,.$$
Note that
$\|\rho\|^2=\|\operatorname{Ric}\|^2$ and
$\|\mathcal{H}_f\|^2=\|\operatorname{Hess}_f\|^2$. 

The triple $(M,g,f)$ is said to be
a \emph{Lorentzian gradient Ricci soliton} if $f$
satisfies the {\it gradient Ricci soliton equation}:
\begin{equation}\label{eq:ricci-soliton}
\operatorname{Hess}_f+\rho=\lambda\, g\text{ for some }\lambda\in\mathbb{R}\,.
\end{equation}
Setting $f=0$ yields the {\it Einstein equation} $\rho=\lambda g$; thus
Equation~(\ref{eq:ricci-soliton}) is a natural generalization of the Einstein equation
and a gradient Ricci soliton can be thought of as a
generalized Einstein manifold. 
Gradient Ricci solitons also correspond to self-similar solutions of the {\it Ricci flow}
$\partial_tg(t)=-2\rho_{g(t)}$. 
For these reasons, gradient Ricci solitons have been extensively
investigated in the literature -- see for example the discussion in
\cite{BGG, Calvaruso-Fino, Cao, Onda} and the references therein.
If $\lambda>0$ (resp. $\lambda=0$ or $\lambda<0)$, then $(M,g,f)$
is said to be \emph{shrinking} (resp. \emph{steady} or \emph{expanding}).
We shall assume for the most part that $(M,g)$ is locally homogeneous. This implies
the scalar curvature is constant. 

One has canonical examples which play a central role in the theory.
Let $(N,g_N)$ be an Einstein manifold with
{\it Einstein constant} $\lambda$, i.e. $\rho_N=\lambda\,g_N$.
Let $M=N\times\mathbb{R}^k$ have the product metric  $g_M$ and let 
$f(x):=\frac{\lambda}{2}\|\pi(x)\|^2$ where $\pi$ is projection on the second factor.
Then $(M,g_M,f)$ is a gradient Ricci soliton and is said to be \emph{rigid}.
Since we are interested in questions of local geometry, by an abuse of notation
we shall also say that $(M,g_M,f)$ is {\it rigid} if $(M,g_M,f)$ is
isomorphic to an open subset of a product $N\times\mathbb{R}^k$ which is rigid.
We shall use the following results of Petersen and Wylie \cite{PW}. Assertion~(2) was first proved in the Riemannian setting but extends easily to arbitrary
signature.

\begin{theorem}\cite{PW}\label{T1}
\ \begin{enumerate}
\item Any locally homogeneous Riemannian gradient Ricci soliton is rigid.
\item Let $(M,g)=(M_1\times M_2,g_1\oplus g_2)$ be 
the direct product of two pseudo-Riemannian manifolds. If $f$ satisfies the gradient Ricci soliton
equation on $(M,g)$, then
$f(x_1+x_2)=f_1(x_1)+f_2(x_2)$ where $f_1$ and $f_2$ satisfy the gradient
Ricci soliton equation on $(M_1,g_1)$ and on $(M_2,g_2)$ separately.
\end{enumerate}\end{theorem}

Assertion (1) was originally proven for homogeneous manifolds, but the assumption of homogeneity can be weakened to local homogeneity by modifying the argument in \cite{PW} Proposition 1 as in the proof of Lemma 2\,(2c). Since  any locally homogeneous Riemannian gradient Ricci soliton is rigid,
the classification is complete in this context. However the possible
geometries are much richer in the Lorentzian setting owing to
the existence of degenerate parallel line fields. For example,
in Example~\ref{Exm1}, we shall present results of \cite{BBGG} showing
that Cahen-Wallach symmetric spaces admit steady non-rigid
gradient Ricci solitons.

\subsection{Outline of the paper and summary of results}
In Section~\ref{S1.2}, we state Lemma~\ref{L2}.
This Lemma, which will be proved in Section~\ref{S2},
summarizes the relevant
results we shall need concerning gradient Ricci solitons with constant
scalar curvature; many of these results rely upon earlier papers. The
analysis there will be local in nature and will rely on the investigation of the
gradient Ricci soliton Equation~(\ref{eq:ricci-soliton}) as this links the
geometry of the manifold, through its Ricci curvature, with the extrinsic geometry
of the level sets of the potential function by means of their second fundamental form.
The signature of the manifold plays no role in Lemma~\ref{L2} and is
completely general. We shall see that if the scalar curvature is constant,
then any solution of \eqref{eq:ricci-soliton} is an
isoparametric function, i.e.
$$\|\nabla f\|^2=b(f)\text{ and }\Delta f=a(f)\text{ for }a,b\text{ smooth on }
\operatorname{Range}(f)\,.$$ 

For the remainder of the paper we shall assume (unless otherwise noted)
that the underlying manifold 
$(M,g)$ is a locally homogeneous Lorentzian manifold and that 
$(M,g,f)$ is a gradient Ricci soliton.
In Section~\ref{S1.3}, we present our results in Theorems~\ref{T3}--\ref{T5} concerning non-steady 
solitons ($\lambda\ne0$);
these results will be proved in Section~\ref{S3}.
In low dimensions, such solitons are rigid; in arbitrary dimensions, 
the eigenvalue structure of the Ricci operator agrees with the corresponding
eigenvalue structure of a rigid soliton, i.e. there are only two eigenvalues $\{0,\lambda\}$. 
In Section~\ref{S1.4}, we present our results concerning steady solitons 
($\lambda=0$)
in Theorems~\ref{T7}--\ref{T8}; 
these will be proved in Section~\ref{S4}. Theorem~\ref{T7} gives
a complete classification if $\|\nabla f\|^2<0$.
In Theorem~\ref{T8}, we shall examine the situation when $\|\nabla f\|^2=0$ and show
the Ricci tensor is either $2$ or $3$ step nilpotent; the metrics in question are 
pure radiation metrics with parallel rays \cite{leistner-nurowski}.
If we further restrict the geometry,
stronger results are available. 
In Section~\ref{S1.5}, we give a complete classification of symmetric
Lorentzian gradient Ricci solitons in Theorem~\ref{T11}. This result will be proved in Section~\ref{S5}.
In Section~\ref{S1.6} in Theorem~\ref{T15}, we give a complete
classification of 3-dimensional Lorentzian 
locally homogeneous gradient Ricci solitons; there are 3 non-trivial
families of examples. Theorem~\ref{T15} will be proved in
Section~\ref{S6}.

The fact that $(M,g)$ is Lorentzian plays a crucial role in many arguments. For
example, when we study the non-steady case, there exists a distinguished null
parallel vector field and there do not exist orthogonal null vector fields -- this is a
Lorentzian phenomena not present in the Riemannian or the higher signature setting.
The fact that $(M,g)$ is locally homogeneous is not simply used to
 ensure that the scalar curvature
is constant, it plays a role in many proofs where we take frame fields consisting at
least in part of Killing vector fields. As our discussion is local in nature, it is not necessary to 
impose global conditions such as global homogeneity or completeness.

\subsection{Consequences of the gradient Ricci soliton equation}\label{S1.2}
Let $\tau$ be the scalar curvature. Let $\nabla f$ be the vector field
which is dual to the $1$-form $df$. It is characterized by the identity
\begin{equation}\label{Eq2}
g(\nabla f,X)=X(f)\text{ for any vector field }X.
\end{equation}
Let $\mathcal{L}$ be the {\it Lie derivative}; 
a vector field $X$ on $(M,g)$ is Killing  if $\mathcal{L}_Xg=0$; $X$ is Killing if and only if
\begin{equation}\label{Eq3}
g(\nabla_XZ,X)=0\text{ for any vector field }X.
\end{equation}
We say $(M,g,f)$ is {\it isotropic} if $\|\nabla f\|^2=0$.
Section~\ref{S2} is devoted to the proof of the following quite general result
concerning gradient Ricci solitons with constant scalar curvature in arbitrary signature.

\begin{lemma}\label{L2}
Let $(M,g,f)$ be a gradient Ricci soliton with constant scalar curvature. 
\begin{enumerate}
\item We have the following relations:
\begin{enumerate}
\item $\operatorname{Ric}(\nabla f)=0$.
\item $\|\nabla f\|^2-2\lambda f=\operatorname{const}$.
\item $R(X,Y,Z,\nabla f)=(\nabla_X \rho)(Y,Z)-(\nabla_Y \rho)(X,Z)$.
\item $\left(\nabla_{\nabla f}\operatorname{Ric}\right)
+\operatorname{Ric}\circ \mathcal{H}_f=R(\nabla f,\cdot)\nabla f$.
\end{enumerate}
\item Let $X$ be a Killing vector field.
\begin{enumerate}
\item $\mathcal{L}_X\left(\operatorname{Hess}_{f}\right)=\operatorname{Hess}_{X(f)}$.
\item  $\grad\{X(f)\}$ is a parallel vector field. 
\item If $\lambda\neq 0$, then $\grad\{X(f)\}=0$ if and only if $X(f)=0$.
\end{enumerate}
\item $\lambda ((n+2)\lambda - \tau)=\|\operatorname{Hess}_f\|^2$.
\item If $(M,g,f)$ is isotropic
and non-steady, then $(M,g)$ is Einstein.
\item If $(M,g,f)$ is steady, then $\|\operatorname{Hess}f\|^2=0$ and $\|\nabla f\|^2=\mu$ is constant.
\end{enumerate}
\end{lemma}

We shall apply different techniques in what follows
 to study the steady and the non-steady cases
since setting  $\lambda\neq 0$ or $\lambda=0$ in Lemma~\ref{L2} gives
significantly different information about the potential function $f$.
By Lemma~\ref{L2}, any isotropic
non-steady gradient Ricci soliton with constant
scalar curvature is Einstein.
However, there exist isotropic steady gradient Ricci solitons 
which are not Einstein \cite{BBGG}. 

\subsection{Non-steady locally homogeneous Lorentzian gradient Ricci solitons}\label{S1.3}

We say that a Lorentzian manifold $(M,g)$
is {\it irreducible} if the holonomy representation has no non-trivial invariant 
subspace and that $(M,g)$ is
{\it indecomposable} if the metric on any non-trivial subspace fixed by the holonomy
representation is degenerate and thus the holonomy representation does not decompose
as a non-trivial direct sum of subrepresentations. The distinction between irreducible and
indecomposable is only relevant in the indefinite setting.
We shall establish the following results in Section~\ref{S3}:
\begin{theorem}\label{T3}
Let $(M,g,f)$ be a locally homogeneous  Lorentzian 
non-steady gradient Ricci soliton. Then one of the following holds: 
\begin{enumerate}
\item $(M,g)$ is irreducible and Einstein.
\item $(M,g,f)$ is rigid, this is, there is a local splitting 
$(M,g,f)=(N\times\mathbb{R}^k_\nu,g_N+g_e,f_N+f_e)$
where $(N,g_N)$ is Einstein with Einstein constant $\lambda$ and $(\mathbb{R}^k_\nu,g_e,f_e)$ is pseudo-Euclidean space, $\nu=0,1$, with $f_e(x):=\frac{\lambda}{2}\|{ x}\|^2$.
\item $(M,g,f)$ locally splits as
\[
(M,g,f)=(N_0\times N_1\times\mathbb{R}^k,g_0+g_1+g_e,f_0+f_1+f_e)
\]
where $(N_0,g_0,f_0)$ is an indecomposable locally homogeneous Lorentzian gradient Ricci soliton, $(N_1,g_1)$ is a Riemannian Einstein manifold with Einstein constant $\lambda$ and $(\mathbb{R}^k,g_e,f_e)$ is Euclidean space with $f_e(x):=\frac{\lambda}{2}\|{ x}\|^2$.
\end{enumerate}
\end{theorem}

We now focus on the situation in Assertion~(3) above and study the indecomposable factor. Recall that a Lorentzian
manifold is said to be {\it Walker} if it admits a parallel null line field,  and {\it strict Walker} if this distribution is spanned by a parallel null vector field;
we refer to \cite{walker-metrics} for further details. We shall say that
$(M,g)$ has {\it harmonic Weyl tensor} if the {\it Schouten tensor} $S$
is {\it Codazzi}. This means (see \cite{besse}):
$$
\nabla_XS_{YZ}=\nabla_YS_{XZ}\text{ where }
S =\rho - \frac{\tau}{2 (n+1)}g\,.
$$

\begin{theorem}\label{T4}
Let $(M,g,f)$ be a locally homogeneous indecomposable Lorentzian non-steady gradient Ricci
soliton which is not Einstein.
\begin{enumerate}
\item Locally, there
exists a Killing vector field $X$ so $U:=\grad\{X(f)\}$ 
is a non-trivial parallel null vector field; thus $(M,g)$ is  strict Walker. 
\item $U$ is unique up to scale, 
$\mathcal{V}:=\{U,\nabla f\}\subset\ker\{\operatorname{Ric}\}$ is a $U$-parallel
Lorentzian distribution,
and $\grad\{U(f)\}=\lambda U$. 
\item $\nabla_U \operatorname{Ric}=\nabla_U\mathcal{H}_f=0$,
$\operatorname{Spec}\{\operatorname{Ric}\}=\operatorname{Spec}\{\mathcal{H}_f\}=\{0,\lambda\}$,
$\operatorname{Ric}$ and $\mathcal{H}_f$ are diagonalizable,  
$ \ker\{\operatorname{Ric}\}=\operatorname{Image}\{\mathcal{H}_f\}$,
and $\ker\{\mathcal{H}_f\}=\operatorname{Image}\{\operatorname{Ric}\}$.
\item The Weyl tensor of $(M,g)$ is harmonic if and only if $(M,g,f)$ is rigid.
\item If $\dim(\ker\{\operatorname{Ric}\})=2$, then $( M,g,f)$ is rigid.
\end{enumerate}
\end{theorem}

This leads to the following classification result in low dimensions:
\begin{theorem}\label{T5}
Let $(M,g,f)$ be a locally homogeneous  Lorentzian 
non-steady gradient Ricci soliton of dimension $m\le4$. Then $(M,g,f)$ is rigid.
\end{theorem}

\begin{remark}
What is indeed proven in Theorem~\ref{T5} is that if the factor $N_0$ of the decomposition given in Theorem~\ref{T3} above is of dimension $n_0\leq 4$ then the gradient Ricci soliton is rigid.
\end{remark}

\subsection{Steady locally homogeneous Lorentzian gradient Ricci solitons}\label{S1.4}
The geometry of the level sets of the potential function plays an essential role in
our analysis; the norm  $\|\nabla f\|^2$ is important as this controls
the nature of the metric on the level sets. The $2$-dimensional case
is trivial; one has \cite{BGG,RicciFlow}:
\begin{theorem}\label{T6}
A steady locally homogeneous Ricci soliton of dimension $2$ either in the
Riemannian or in the Lorentzian setting is flat.
\end{theorem}

The following two results will be established in Section~\ref{S4}:

\begin{theorem}\label{T7}
Let $(M,g,f)$ be a locally homogeneous steady gradient Lorentzian Ricci soliton.
If $\|\nabla f\|^2<0$, then
$(M,g)$ splits locally as an isometric
product $(\mathbb{R}\times N,-dt^2+g_N)$, 
where $(N,g_N)$ is a flat Riemannian manifold and $f$ is orthogonal projection on 
$\mathbb{R}$.
\end{theorem}

The cases when $\|\nabla f\|^2\geq 0$ are less rigid in the steady setting.
Several examples in the spacelike case $\|\nabla f\|^2>0$ are known \cite{BBGG,BGG},
but little more of a general nature is known about this case.
In the isotropic case one has some restrictions on the Ricci operator;
in particular, it must be nilpotent. Recall that a tensor $T$ is said to be {\it recurrent}
if there is a smooth $1$-form $\omega$ so that $\nabla_XT=\omega(X)T$.

\begin{theorem}\label{T8}
Let $(M,g,f)$ be an isotropic locally homogeneous  Lorentzian steady gradient Ricci soliton. One
of the following two possibilities pertains:
\begin{enumerate}
\item  $\mathcal{H}_f=-\operatorname{Ric}$ has rank $2$ and is $3$-step nilpotent.
\item $\mathcal{H}_f=-\operatorname{Ric}$ has rank $1$ and is $2$-step nilpotent. 
In this case $(M,g)$ is locally a strict Walker manifold, more specifically:
\begin{enumerate}
\item $\ker\{\mathcal{H}_f\}=\nabla f^\perp$ and 
$\operatorname{Image}\{\mathcal{H}_f\}=\nabla f$.
\item $\nabla f$ is a recurrent vector field and $\nabla f^\perp$ is an integrable
totally geodesic distribution with leaves the level sets of $f$.
\item Let $P\in M$. At least one of the following possibilities holds near $P$.
\begin{enumerate} 
\item There exists a Killing vector field $F$ so $\grad\{F(f)\}$
is a null parallel vector field.
\item There exists a smooth function $\psi$ defined near $P$ so\
$\psi\ \nabla f$ is a null parallel vector field.
\end{enumerate}
\end{enumerate}\end{enumerate}
\end{theorem}

We shall illustrate possibility~(2) in Example~\ref{Exm1} presently.

\subsection{Symmetric Lorentzian gradient Ricci solitons}\label{S1.5}
Stronger results are available if $(M,g)$ is {\it locally symmetric}; this implies $\nabla R=0$. 
\begin{definition}\label{D9}
We say that $(N,g_N)$ is a {\it Cahen-Wallach symmetric space} if
there are coordinates $(t,y,x_1,\dots,x_n)$ so:
\begin{equation}\label{E4}
g=2\, dt dy+\left(\sum_{i=1}^{n} \kappa_i\, x_i^2\right) dy^2+\sum_{i=1}^{n} dx_i^2
\text{ for }0\ne\kappa_i\in\mathbb{R}\,.
\end{equation}
We shall always assume that all $\kappa_i\ne0$ to ensure that $(N,g_N)$ is indecomposable.
\end{definition}

We refer to \cite{CW70,CLPTV1} for the proof of Assertion~(1) in the following
result and to \cite{BBGG} for the proof of Assertion~(2) in the following result:
\begin{theorem}\label{T10}\
\begin{enumerate}
\item Let $(M,g)$ be a Lorentzian locally symmetric space.
\begin{enumerate}
\item If $(M,g)$ is irreducible, then $(M,g)$ has constant sectional curvature.
\item If $(M,g)$ is indecomposable but reducible, then $(M,g)$
is a Cahen-Wallach symmetric space.\end{enumerate}
\item If $(M,g,f)$ is a Cahen-Wallach 
gradient Ricci soliton, then
$(M,g,f)$ is steady, $f=a_0+a_1y+\frac14\sum_i\kappa_iy^2$, and $\nabla f=(a_0+\frac12\sum_i\kappa_iy)\partial_t$
is null.
\end{enumerate}
\end{theorem}

Theorem~\ref{T10} will play a crucial
role in the proof that we shall give of the following result in Section~\ref{S5}:

\begin{theorem}\label{T11}
Let $(M,g,f)$ be a locally symmetric Lorentzian gradient Ricci soliton. Then $(M,g)$ splits locally as a product $M=N\times\mathbb{R}^k$ where
\begin{enumerate}
\item
if $(M,g,f)$ is not steady, then $(N,g_N)$ is Einstein and the soliton is rigid,
\item
if $(M,g,f)$ is steady, then $(N,g_N,f_N)$ is locally isometric to a Cahen-Wallach symmetric space.
\end{enumerate}
\end{theorem}

\subsection{Three-dimensional locally homogeneous gradient Ricci solitons}\label{S1.6}
We will establish the following 2 results in $3$-dimensional geometry in Section~\ref{S6}.
Let $(M,g)$ be a Lorentzian manifold of dimension $3$. We suppose first that $(M,g)$ is strict
Walker, i.e. admits a null parallel vector field. We may then 
(see, for example, \cite{walker-metrics}) find local adapted coordinates
$(t,x,y)$ so that
\begin{equation}\label{E5x}
g= 2 dt dy +  dx^2+\phi(x,y)dy^2.
\end{equation}
The following is of independent interest; we drop for the moment the assumption that
the metric is locally homogeneous and focus on Walker geometry:

\begin{theorem}\label{T12}
Let $(M,g)$ be a non-flat 3-dimensional Lorentzian strict Walker manifold. 
Then $(M,g,f)$ is a gradient Ricci soliton if and only if there exist a cover of $M$ by coordinate
systems where the metric has the form given in Equation~(\ref{E5x}) where one of the following occurs:
\begin{enumerate}
\item $\phi(x,y)=\frac{1}{\alpha^2}\,a(y)\,e^{\alpha x}+x\, b(y)+c(y)$ and
$f(x,y)=x\, \alpha+\gamma(y)$ where $\alpha\in\mathbb{R}$ and
 $\gamma''(y)=-\frac{1}{2}\, \alpha b(y)$. In
this setting, $\nabla f=\alpha\partial_x+\gamma^\prime(y)\partial_t$ is spacelike.
\item $\phi(x,y)=x^2\, a(y)+x\, b(y)+ c(y)$ and $f(x,y)=\gamma(y)$ where
$\gamma''(y)=\frac{1}{4}\, a(y)$.
In this setting $\nabla f=\gamma^\prime\partial_t$ is null.
\end{enumerate}
Moreover, in both cases the Ricci soliton is steady.
\end{theorem}

\begin{definition}\label{D13}\rm
Adopt the notation of Equation~(\ref{E5x}).
\begin{enumerate}
\item Let $\phi(x,y)=b^{-2}e^{bx}$ for $0\ne b\in\mathbb{R}$ define $\mathcal{N}_b$.
\item Let $\phi(x,y)=\frac12x^2\alpha(y)$
where $\alpha_y(y)=c\alpha^{3/2}(y)$ and $\alpha(y)>0$ define $\mathcal{P}_c$ .
\item Let $\phi(x,y)=\pm x^2$ define the Cahen-Wallach symmetric space $\mathcal{CW}_\pm$. 
\end{enumerate}
\end{definition}

The following result was established in \cite{GR-G-N}:
\begin{theorem}\label{T14}
Let $(M,g)$ be a locally homogeneous Lorentzian strict Walker  manifold of dimension $3$. Then $(M,g)$
is locally isometric to one of the manifolds given in {\rm Definition~\ref{D13}}.
\end{theorem}

We can now state our classification result:

\begin{theorem}\label{T15}
Let $(M,g,f)$ be a Lorentzian locally homogeneous gradient Ricci soliton
of dimension $3$. If $(M,g,f)$ is 
non-trivial, then either it is rigid or $(M,g)$ is locally isometric to either $\mathcal{CW}_\pm$, 
$\mathcal{P}_c$ or  $\mathcal{N}_b$ as defined above and the soliton is steady.
Moreover $\nabla f$ is null if $(M,g)=\mathcal{P}_c$ or if $(M,g)=\mathcal{CW}_\pm$, and
$\nabla f$ is spacelike if $(M,g)=\mathcal{N}_b$.
\end{theorem}

\section{ Consequences of the gradient Ricci soliton equation\\
The proof of Lemma~\ref{L2}}\label{S2}
In Section~\ref{S2.1}, we establish Assertion~(1), in Section~\ref{S2.2} we derive Assertion~(2),
in Section~\ref{S2.3}, we prove Assertion~(3), in Section~\ref{S2.4}, we verify Assertion~(4),
and in Section~\ref{S2.5}, we complete the proof of Lemma~\ref{L2} by checking Assertion~(5).

\subsection{The proof of Lemma~\ref{L2}~(1)}\label{S2.1}
If $(M,g,f)$ is a gradient Ricci soliton, then  
$\nabla \tau= 2\operatorname{Ric}(\nabla f)$ \cite{E-LN-M,PW2}. 
Assertion~(1a) now follows as $\nabla\tau=0$. We also have
\cite{BGG,Cao,E-LN-M,PW2} that
$\tau+\|\nabla f\|^2-2\lambda f=\operatorname{const}$; Assertion~(1b) now follows. We
refer to \cite{BGG,FL-GR} for the proof of Assertion~(1c) which holds without assuming
$\tau=\operatorname{const}$.
The identity
$$
\left(\nabla_{\nabla f}\operatorname{Ric}\right)+\operatorname{Ric}\circ\mathcal{H}_f
=R(\nabla f,\cdot)\nabla f+\textstyle\frac12 \nabla \nabla\tau
$$
was proved in the Riemannian setting \cite{PW2}. One can
use analytic continuation to extend this identity to the indefinite setting (or simply observe the proof
goes through without change in the higher signature context).
Assertion~(1d) now follows once again using the
fact that $\tau$ is constant.\hfill
\qed

\subsection{The proof of Lemma~\ref{L2}~(2)}\label{S2.2}
Let $X$ be a Killing vector field.
Fix a point $P$ of $M$ so that $X(P)\ne0$; 
Assertion~(2) for $P$ where $X(P)=0$ will then follow by continuity. Choose a system
of local coordinates $(x_1,\dots,x_{n+2})$ so that $X=\partial_{x_1}$. Set
$g_{ij}:=g(\partial_{x_i},\partial_{x_j})$ and observe that
\begin{eqnarray*}
\partial_{x_1}\, g_{ij}&=&g(\nabla_{\partial_{x_1}}\partial_{x_i},\partial_{x_j})+g(\partial_{x_i},\nabla_{\partial_{x_1}}\partial_{x_j})\\&=&g(\nabla_{\partial_{x_i}}\partial_{x_1},\partial_{x_j})+g(\partial_{x_i},\nabla_{\partial_{x_j}}\partial_{x_1})
=(\mathcal{L}_{\partial_{x_1}}\, g)(\partial_{x_i},\partial_{x_j})\,.
\end{eqnarray*}
Thus $\partial_{x_1}\, g_{ij}=0$ so
$\partial_{x_1}\, \Gamma_{ij}{}^k=0$ as well.
We establish Assertion~(2a) by computing:
\begin{eqnarray*}
&&(\mathcal{L}_{\partial_{x_1}}\operatorname{Hess}_f)(\partial_{x_i},\partial_{x_j})
=\mathcal{L}_{\partial_{x_1}}\operatorname{Hess}_f(\partial_{x_i},\partial_{x_j})\\
&=&\mathcal{L}_{\partial_{x_1}}\left(\partial^2_{x_ix_j}(f)-\Gamma_{ij}{}^k\partial_{x_k}(f)\right)\\
&=&\partial^3_{x_1x_ix_j}(f)-\partial_{x_1}(\Gamma_{ij}{}^k)\partial_{x_k}(f)-\Gamma_{ij}{}^k\partial^2_{x_1x_k}(f)\\
&=&\partial^2_{x_ix_j}\partial_{x_1}(f)-\Gamma_{ij}{}^k\partial_{x_k}\partial_{x_1}(f)
=\operatorname{Hess}_{\partial_{x_1}(f)}(\partial_{x_i},\partial_{x_j})\,.
\end{eqnarray*}
Since $\mathcal{L}_Xg=0$ and since $\rho$ is natural,
$\mathcal{L}_X\rho=0$. Equation \eqref{eq:ricci-soliton} implies that
 $\mathcal{L}_X\operatorname{Hess}_{f}=0$, and therefore by Assertion~(2a),
 $\operatorname{Hess}_{X(f)}=0$. 
 Consequently, $\grad\{X(f)\}$ is parallel. This establishes Assertion~(2b).
Assume now that $\lambda\neq 0$. It is clear that $\grad\{X(f)\}=0$ if $X(f)=0$. 
Conversely, if $\grad\{X(f)\}=0$, then $X(f)=\kappa$ for some constant $\kappa$. 
Since the scalar curvature is constant, Assertion~(1) implies
that $\operatorname{Ric}(\nabla f)=0$. Since $X$ is a Killing vector field,
\begin{eqnarray*}
0&=&\nabla f(\kappa)=\nabla f(X(f))= \nabla f\,g(\nabla f,X)
=g(\nabla_{\nabla f}\nabla f,X)+g(\nabla f,\nabla_{\nabla f}X)\\
&=&\operatorname{Hess}_f(\nabla f,X)+\frac{1}{2}(\mathcal{L}_Xg)(\nabla f,\nabla f)
=-\rho(\nabla f,X)+\lambda\, g(\nabla f,X)=\lambda\, \kappa\,.
\end{eqnarray*}
Thus $\kappa=0$. Consequently $\grad\{X(f)\}=0$ if and only if $X(f)=0$.
This establishes Assertion~(2c).\hfill
\qed

\subsection{The proof of Lemma~\ref{L2}~(3)}\label{S2.3}
We have the Bochner identity:
\begin{equation}\label{E6}
{\textstyle\frac12}\,\Delta\, g(\nabla f,\nabla f)
=\|\operatorname{Hess}_f\|^2 + \rho(\nabla f,\nabla f) + g(\nabla\Delta f,\nabla f).
\end{equation}
By Assertion~(1), $\operatorname{Ric}(\nabla f)=0$ and
$\|\nabla f\|^2-2\lambda f=\operatorname{const}$. Thus the left-hand side of
Equation~(\ref{E6}) becomes
$\frac{1}{2}\,\Delta\, g(\nabla f,\nabla f)
=\lambda\, \Delta f-\frac{1}{2}\,\Delta\tau$. Taking
the trace in Equation \eqref{eq:ricci-soliton} shows that 
$\Delta f=(n+2)\lambda -\tau$ and hence
$\frac{1}{2}\,\Delta\, g(\nabla f,\nabla f)=\lambda((n+2)\lambda - \tau)$. 
On the other hand, since $\operatorname{Ric}(\nabla f)=0$ and 
$\nabla\Delta f=-\nabla\tau=0$, the right-hand side in Bochner formula reduces to  
$\|\operatorname{Hess}_f\|^2$.\qed

\subsection{The proof of Lemma~\ref{L2}~(4)}\label{S2.4}
If $\|\nabla f\|^2=0$, we may apply Assertion~(1) to see 
$2\lambda f=\operatorname{const}$. 
Since $\lambda\neq 0$, $f$ is constant and $(M,g)$ is Einstein. \qed

\subsection{The proof of Lemma~\ref{L2}~(5)}\label{S2.5}
If $\lambda=0$, then
 $\|\operatorname{Hess}_f\|^2=0$. By Equation~(\ref{eq:ricci-soliton}),
$\mathcal{H}_f=-\operatorname{Ric}$ and thus
$\operatorname{Ric}(\nabla f)=0$ implies $\mathcal{H}_f(\nabla f)=0$. Consequently
$\nabla f$ is a geodesic vector field.
Next, using the identity $\tau+\|\nabla f\|^2-2\lambda f=\operatorname{const}$, 
one has that $\|\nabla f\|^2$ is constant and therefore $f$ is a solution of the 
{\it Eikonal equation $\|\nabla f\|^2=\mu$}.
\qed

\section{Non-steady locally homogeneous gradient Ricci solitons\\
the proof of Theorems~\ref{T3}--\ref{T5}}\label{S3}

By Lemma~\ref{L2}, isotropic non-steady locally homogeneous gradient Ricci solitons
are Einstein. Consequently, we shall concentrate henceforth on the study of 
non-isotropic non-steady locally homogeneous gradient Ricci solitons.
In Section~\ref{S3.1}, we will prove Theorem~\ref{T3}, in Section~\ref{S3.2}, we will establish
Theorem~\ref{T4}, and in Section~\ref{S3.3}, we will establish Theorem~\ref{T5}.
We shall use Lemma~\ref{L2} repeatedly and without further reference in what follows.
Throughout Section~\ref{S3}, we shall let $(M,g,f)$ be a locally homogeneous non-steady
gradient Ricci soliton.

\subsection{The proof of Theorem~\ref{T3}}\label{S3.1} Assume that $(M,g)$
is irreducible or, equivalently, that there are no non-trivial parallel distributions on $M$.
Consequently any parallel vector field is trivial. Let $X$ be a Killing vector
field. Then $\grad\{X(f)\}$ is a parallel vector field and thus
$\grad\{X(f)\}=0$ so $X(f)$ is constant  and hence $X(f)=0$. 
Since the underlying Lorentzian structure $(M,g)$ is locally homogeneous,
there are $(n+2)$ linearly independent Killing vector fields
$X_1,\dots,X_{n+2}$ locally. 
Consequently $f$ is constant and the metric is Einstein.
This establishes Assertion~(1) of Theorem~\ref{T3}.

We now apply the local splitting result of Assertion~(2) in Theorem~\ref{T1}. 
Let $X$ be a Killing vector field on $(M,g)$. 
If $\grad\{X(f)\}$ is spacelike or timelike, then we may split,
at least locally, a one-dimensional factor 
from $(M,g)$ and decompose locally
$$(M,g,f)=(N\times\mathbb{R},g_N\oplus g_e,f_N+ f_e)\,.$$
If $\grad\{X(f)\}$ is timelike, then $(N,g_N)$ is Riemannian and by Assertion~(1) of Theorem~\ref{T1},
rigid which would finish the discussion.
Thus we may assume $(N,g_N)$ is Lorentzian so $\grad\{X(f)\}$ is spacelike and the
factor $(\mathbb{R},g_e)$ is
positive definite. We proceed inductively to decompose
$(M,g,f)=(N\times\mathbb{R}^k,g_N\oplus g_e,f_N+ f_e)$ (at least locally) so that 
$(N,g_N,f_N)$ is a locally homogeneous Lorentzian Ricci soliton with $\grad\{X(f)\}$ null or zero for all Killing vector fields $X$. Now two possibilities may occur. If $N$ is indecomposable, Assertion~(3) follows with trivial $N_1$. If $N$ is decomposable, then either $N$ is Einstein  and Assertion~(2) holds (this is the case if $\grad\{X(f)\}=0$ for all Killing vector fields in $N$) or $N$ decomposes as $N=N_0\times N_1$ where $N_0$ is Lorentzian and indecomposable (the latter happens if there exists a Killing vector field $X$ so that $\grad\{X(f)\}$ is null). $(N_1,g_1,f_1)$ is a Riemannian locally homogeneous gradient Ricci soliton which, as a consequence of Theorem~\ref{T1}, is Einstein. This establishes Theorem~\ref{T3}.\hfill\qed

\subsection{The proof of Theorem~\ref{T4}}
\label{S3.2}
We establish Assertions~(1)-(5) of Theorem~\ref{T4} seriatim. We suppose
$(M,g)$ is not decomposable and is not Einstein. 
\subsubsection*{The proof of Theorem~\ref{T4}~(1)}
We must show  there exists $X$ so $U=\grad\{X(f)\}$ is a parallel null vector field.
Let $Z$ be any Killing vector field.
Since $(M,g)$ is not decomposable and since $\grad\{Z(f)\}$ is parallel, 
$\grad\{Z(f)\}$ must be isotropic. If $\grad\{Z(f)\}$ vanishes
for all such $Z$, then $f$ is constant and hence $(M,g)$ is Einstein which is contrary to our
assumption. Thus $U:=\grad\{Z(f)\}$ has the desired properties for some Killing vector field $Z$.

\subsubsection*{The proof of Theorem~\ref{T4}~(2)} We must show that $U$ is
unique up to scale, that $U\in\ker\{\operatorname{Ric}\}$, and that $\grad\{U(f)\}=\lambda U$.
Suppose that there are
two Killing vector fields $Z_1$ and $Z_2$ on $(M,g)$ so that $\grad\{Z_1(f)\}$ and $\grad\{Z_2(f)\}$
are linearly independent. Since the signature is Lorentzian, 
$\operatorname{Span}\{\grad\{Z_1(f)\},\grad\{Z_2(f)\}\}$ can not be a null distribution.
Consequently, there exists a linear combination
$Z=a_1Z_1+a_2Z_2$ so $\grad\{Z(f)\}$ is either timelike or spacelike. This implies
that $(M,g)$ is decomposable which is false. Thus the vector field $U=\grad\{Z(f)\}$ is
unique up to scale. 

Since $U$ is parallel, it is Killing and hence $\grad\{U(f)\}=\alpha U$ for 
some $\alpha\in\mathbb{R}$.
We must now show that $\operatorname{Ric}(U)=0$.
Let $\{Z_1,Z_2,\dots,Z_{n+2}\}$ be a local basis of Killing vector fields.
Choose the notation so $Z=Z_1$. We  then have $\grad\{Z_i(f)\}=\mu_iU$ for $i\ge2$. 
Since $\grad\{Z_i(f)\}$ is parallel, necessarily $\mu_i$ is constant. 
By replacing $Z_i$ by $Z_i-\mu_iZ_1$,
we may assume therefore that $\grad\{Z_i(f)\}=0$ for $i\ge2$. 
Since $\lambda\ne0$, Lemma~\ref{L2}
implies $Z_i(f)=0$ for $i\ge2$. We use  Equation~(\ref{eq:ricci-soliton}) and
Equation~(\ref{Eq2}) to see:
\begin{eqnarray}
&&g(U,\nabla f)=g(\operatorname{grad}\{Z_1(f)\},\nabla f)
=g(\operatorname{grad}\{g(Z_1,\nabla f)\},\nabla f)\\
&&\quad=\nabla f\, g(Z_1 ,\nabla f)\nonumber
=g(\nabla_{\nabla f}Z_1,\nabla f)+g(Z_1,\nabla_{\nabla f}\nabla f)
\label{Ex5}\\
&&\quad=\operatorname{Hess}_f(Z_1,\nabla f)
=\lambda g(Z_1,\nabla f)=\lambda Z_1(f)\neq 0\,,\nonumber
\end{eqnarray}
where by Equation~(\ref{Eq3}),
$g(\nabla_{\nabla f}Z_1,\nabla f)=0$ since $Z_1$ is Killing.
As $g(U,\nabla f)\ne0$ and as $U$ is a null vector, 
$\mathcal{V}:=\operatorname{Span}\{U,\nabla f\}$ has Lorentzian signature.
We have that $\grad\{U(f)\}\neq 0$ due to Lemma~\ref{L2} so  $\alpha\ne0$.

If $X$ is an arbitrary vector field, we study $\mathcal{H}_f(U)$ 
by computing:
\[
\operatorname{Hess}_f(X,U)=g(U, \nabla_{X} \nabla f)
=X g(U, \nabla f)=g(X,\grad\{U(f)\})
=\alpha g(X,U)\,.
\]
This shows that $\mathcal{H}_f(U)=\alpha U$. 
Since $\mathcal{H}_f(\nabla f)=\lambda \nabla f$, we also have:
\[
\alpha g(\nabla f,U)=\operatorname{Hess}_f(\nabla f,U)=\lambda g(\nabla f,U)
\]
so $\alpha=\lambda$. 
By Equation~\eqref{eq:ricci-soliton}, $\operatorname{Ric}(U)=0$.
Since $\nabla_U U=0$ and $\nabla_U \nabla f=\lambda U$, $\nabla_U$ 
preserves $\mathcal{V}\subset\ker\{\operatorname{Ric}\}$. 
This proves Assertion~(2).

\subsubsection*{The proof of Theorem~\ref{T4}~(3)}
We have shown that
$\mathcal{V}:=\operatorname{Span}\{U,V\}\subset
\ker\{\operatorname{Ric}\}$ is a $U$-parallel
Lorentzian distribution. Consequently
$\mathcal{V}^\perp$ is a $\operatorname{Ric}$ invariant distribution with a positive definite signature.
Since $\operatorname{Ric}$ is self-adjoint, 
there exists an orthonormal basis $\{E_1,\dots, E_n\}$ of 
$\mathcal{V}^\perp$ so $\operatorname{Ric}(E_i)=\alpha_i E_i$; 
the $\alpha_i$ are constant
since $(M,g)$ is locally homogeneous. This proves in particular that 
$\operatorname{Ric}$ and
$\mathcal{H}_f=\lambda\,\operatorname{Id}-\operatorname{Ric}$ are diagonalizable. We now show
 that $\nabla_U$ preserves the eigenspaces in $\mathcal{V}^\perp$.
For $i\neq j$, since $U$ is parallel $R(U,E_i,E_j,\nabla f)=0$. 
By Lemma~\ref{L2}~(1):
\begin{eqnarray*}
0&=&R(U,E_i,E_j,\nabla f)= (\nabla_U \rho)(E_i,E_j)-(\nabla_{E_i} \rho)(U, E_j)\\
&=& U \rho(E_i,E_j)-\rho(\nabla_U E_i,E_j)-\rho(E_i,\nabla_U E_j)\\
&&\quad-E_i \rho(U,E_j)+\rho(\nabla_{E_i} U, E_j)+\rho(U, \nabla_{E_i} E_j)\\
&=& -\alpha_j g(\nabla_U E_i,E_j)-\alpha_i g(E_i,\nabla_U E_j)\\
&=&(\alpha_i -\alpha_j) g(\nabla_UE_i,E_j)\,.
\end{eqnarray*}
We conclude that if $E_i$ and $E_j$ belong to different eigenspaces $\nabla_U E_i$ 
is orthogonal to $E_j$. Hence, $\nabla_U$ commutes with $\operatorname{Ric}$ 
and, as a consequence of the Ricci soliton equation \eqref{eq:ricci-soliton}, 
it also commutes with $\mathcal{H}_f$. Consequently, as desired, 
 $\nabla_U\operatorname{Ric}=0$ and $\nabla_U\mathcal{H}_f=0$.

We must show that $0$ and $\lambda$ 
are the only eigenvalues of $\operatorname{Ric}$.
Normalize $V$ to be a multiple of $\nabla f$ so $g(V,V)=\epsilon=\pm1$. 
Let $S$ be any level set of $f$. The integral curves of $U$ are transversal to $S$ 
because $g(U,\nabla f)\neq 0$. 
Use parallel transport along the integral curves of $U$ to extend
the local frame   $\{E_1,\dots, E_n\}$ from $S$ to a neighborhood of $S$ to define a
local frame field $\{F_1,\dots, F_n\}$ for $\mathcal{V}^\perp$
 such that $\nabla_U F_i =0$.
Since $\nabla_U\operatorname{Ric}=0$,
the vector fields $F_i$ are still eigenvectors of the Ricci operator 
$\operatorname{Ric}$. We shall
use this local frame field
to see that $\operatorname{Ric}$ has only two eigenvalues 
$\{0,\lambda\}$. First note that
\begin{eqnarray*}
&&(\nabla_{\nabla f} \rho)(F_i,F_i)
=\nabla f\,\rho(F_i,F_i)-2\rho(\nabla_{\nabla f}F_i,F_i)\\
&=&\alpha_i \nabla f\,g(F_i,F_i)-2\alpha_i g(\nabla_{\nabla f}F_i,F_i)
= \alpha_i \left(\nabla_{\nabla f} g\right)(F_i,F_i)=0\,.
\end{eqnarray*}
We use Lemma~\ref{L2} to compute:
\medbreak\quad
$\rho(F_i,F_i)= \epsilon R(F_i, V, F_i, V)$
\smallbreak\qquad\quad
$+\sum_{j\neq i} R(F_i,F_j,F_i,V)g(F_j,V)+\sum_{j\neq i} R(F_i,F_j,F_i,F_j)$
\medbreak\qquad
$= \frac{\epsilon}{\|\nabla f\|^2} \left( (\nabla_{F_i} \rho)(\nabla f,F_i)
-(\nabla_{\nabla f} \rho)(F_i,F_i)\right)$
\smallbreak\qquad\quad
$+\sum_{j\neq i} R(F_i,F_j,F_i,V)g(F_j,V)+\sum_{j\neq i} R(F_i,F_j,F_i,F_j)$
\medbreak\qquad
$=\frac{\epsilon}{\|\nabla f\|^2} \left(F_i\rho(\nabla f,F_i)-\rho(\nabla_{F_i} \nabla f,F_i)
-\rho(\nabla f,\nabla_{F_i}F_i)\right)$
\smallbreak\qquad\quad
$+\sum_{j\neq i} R(F_i,F_j,F_i,V)g(F_j,V)+\sum_{j\neq i} R(F_i,F_j,F_i,F_j)$
\medbreak\qquad
$= -\frac{\epsilon}{\|\nabla f\|^2} \rho(\mathcal{H}_f{F_i},F_i)$
\smallbreak\qquad\quad
$+\sum_{j\neq i} R(F_i,F_j,F_i,V)g(F_j,V)+\sum_{j\neq i} R(F_i,F_j,F_i,F_j)$.
\medbreak\noindent
Since we have shown that $\nabla_U \rho=0$, 
we have that $U\rho(F_i,F_i)=2\rho(\nabla_U F_i, F_i)$
which vanishes. 
We now differentiate  the three summands in the previous expression
with respect to $U$:
\medbreak\quad
$U(-\frac{1}{\|\nabla f\|^2} \rho(\mathcal{H}_f{F_i},F_i))
= \frac{ U g(\nabla f,\nabla f)}{\|\nabla f\|^4} \rho(\mathcal{H}_f{F_i} ,F_i)
-\frac{1}{\|\nabla f\|} U\rho(\mathcal{H}_f{F_i},F_i)$
\smallbreak\qquad
$= \frac{2\lambda g(U,\nabla f)}{\|\nabla f\|^4} \rho(\mathcal{H}_f{F_i} ,F_i)$
$-\frac{1}{\|\nabla f\|} \left(\rho(\nabla_U \mathcal{H}_f{F_i},F_i)
+\rho(\nabla_{F_i} \nabla f,\nabla_U F_i)\right)$
\smallbreak\qquad
$=\frac{2\lambda g(U,\nabla f)}{\|\nabla f\|^4} \rho(\mathcal{H}_f{F_i} ,F_i)$
$-\frac{1}{\|\nabla f\|} \left(\rho(\mathcal{H}_f(\nabla_U F_i),F_i)
+\rho(\nabla_{F_i} \nabla f,\nabla_U F_i)\right)$
\smallbreak\qquad
$=\frac{2\lambda g(U,\nabla f)}{\|\nabla f\|^4} \alpha_i(\lambda-\alpha_i)$.
\medbreak\noindent

\medskip

\medbreak\quad
$U\left(R(F_i,F_j,F_i,\nabla f)g(F_j,\nabla f)\right)$
\smallbreak\qquad
$=\big\{\left(\nabla_U\,R\right)(F_i,F_j,F_i,\nabla f)
+R(\nabla_U F_i,F_j,F_i,\nabla f)+R(F_i,\nabla_UF_j,F_i,\nabla f)$
\smallbreak\qquad\quad
$+R(F_i,F_j,\nabla_UF_i,\nabla f)+R(F_i,F_j,F_i,\nabla_U\nabla f)\big\}
g(F_j,\nabla f)$
\smallbreak\qquad\quad
$+ R(F_i,F_j,F_i,\nabla f)\left(g(\nabla_UF_j,\nabla f)
+g(F_j,\nabla_U \nabla f)\right)$
\medbreak\qquad
$=\big\{-\left(\nabla_{F_i}\,R\right)(F_j,U,F_i,\nabla f)
-\left(\nabla_{F_j}\,R\right)(U,F_i,F_i,\nabla f)$
\smallbreak\qquad\quad
$+R(\nabla_U F_i,F_j,F_i,\nabla f)$
$+R(F_i,\nabla_UF_j,F_i,\nabla f)+R(F_i,F_j,\nabla_UF_i,\nabla f)$
\smallbreak\qquad\quad
$+R(F_i,F_j,F_i,\lambda U)\big\}g(F_j,\nabla f)$
\smallbreak\qquad\quad
$+ R(F_i,F_j,F_i,\nabla f)\left(g(\nabla_UF_j,\nabla f)+\lambda g(F_j,U )\right)$
\smallbreak\qquad
$= \big\{R(\nabla_U F_i,F_j,F_i,\nabla f)+R(F_i,\nabla_UF_j,F_i,\nabla f)$
\smallbreak\qquad\quad
$+R(F_i,F_j,\nabla_UF_i,\nabla f)\big\}g(F_j,\nabla f)$
$+R(F_i,F_j,F_i,\nabla f)g(\nabla_UF_j,\nabla f)$
\medbreak\qquad$=0$.
\medbreak\noindent Consequently along the slice $S$ we have:
\medbreak\quad
$U\left(R(F_i,F_j,F_i,V)g(F_j,V)\right)=
U\|\nabla f\|^{-2}R(F_i,F_j,F_i,\nabla f)g(F_j,\nabla f)$
\smallbreak\qquad\quad
$+\|\nabla f\|^{-2}U\left(R(F_i,F_j,F_i,\nabla f)g(F_j,\nabla f)\right)=0$,
\medbreak\quad
$U R(F_i,F_j,F_i,F_j)
=(\nabla_U R)(F_i,F_j,F_i,F_j)+2R(\nabla_U F_i,F_j,F_i,F_j)$
\smallbreak\qquad\quad
$+2 R(F_i,\nabla_U F_j,F_i,F_j)$
\smallbreak\qquad
$=-(\nabla_{F_i} R)(F_j,U,F_i,F_j)-(\nabla_{F_j} R)(U,F_i,F_i,F_j)$
\smallbreak\qquad\quad
$+2R(\nabla_U F_i,F_j,F_i,F_j)+2 R(F_i,\nabla_U F_j,F_i,F_j)$
\smallbreak\qquad
$=2R(\nabla_U F_i,F_j,F_i,F_j)+2 R(F_i,\nabla_U F_j,F_i,F_j)$
\smallbreak\qquad
$=\,0$.

Hence, the following equation holds: 
\[0={2\lambda g(U,\nabla f)}\|\nabla f\|^{-4} \alpha_i (\lambda-\alpha_i).\]
\medbreak\noindent
Since $\lambda$ and $g(U,\nabla f)$ are different from $0$, either 
$\alpha_i=0$ or $\alpha_i=\lambda$ for $i=1,\dots, n$.  
Since the level set $S$ of $f$ which was
chosen was arbitrary, this is true on all of $M$.
By Equation~(\ref{eq:ricci-soliton}) we have
 $\mathcal{H}_f+\operatorname{Ric}=\lambda\,\operatorname{Id}$. The remaining
conclusions of Assertion~(3) are now immediate from the discussion above.
 
\subsubsection{The proof of Theorem~\ref{T4}~(4)}
Recall that $(M,g)$ has  a harmonic Weyl tensor if its Schouten tensor 
$S = \rho - \frac{\tau}{2 (n+1)}g$ is Codazzi, i.e., $\nabla_XS_{YZ}=\nabla_YS_{XZ}$
(see \cite{besse}). If the Weyl tensor is harmonic then 
$(\nabla_X \rho)(Y,Z)-(\nabla_Y \rho)(X,Z)=0$
since the scalar curvature is constant. 
Choose $E_1,E_2\in \operatorname{Image}\{\mathcal{H}_f\}$ and 
$F\in \operatorname{Image}\{\operatorname{Ric}\}$. 
We use Assertion~(3) to compute
\[
0=(\nabla_{E_1} \rho)(F,E_2)-(\nabla_F \rho)(E_1,E_2)
=\rho(F,\nabla_{E_1}E_2)=\lambda g(F,\nabla_{E_1}E_2)\,.
\]
Choose $E\in \operatorname{Image}\{\mathcal{H}_f\}$ and 
$F_1,F_2\in \operatorname{Image}\{\operatorname{Ric}\}$. 
We show the two eigenspaces are parallel and that the soliton
is rigid by computing:
\begin{eqnarray*}
0&=&(\nabla_{F_1} \rho)(E,F_2)-(\nabla_E \rho)(F_1,F_2)\\
&=&\rho(\nabla_{F_1}E,F_2)-E\rho(F_1,F_2)
+\rho(\nabla_EF_1,F_2)+\rho(F_1,\nabla_EF_2)\\
&=& \lambda g(\nabla_{F_1}E,F_2)-\lambda E\,g(F_1,F_2)
+\lambda g(\nabla_EF_1,F_2)+\lambda g(F_1,\nabla_EF_2)\\
&=&\lambda g(\nabla_{F_1}E,F_2)\,.
\end{eqnarray*}

\subsubsection*{The proof of Theorem~\ref{T4}~(5)}
We apply Theorem~\ref{T3}.
If $\dim(\ker\{\operatorname{Ric}\})=2$, then
$\mathcal{V}=\ker\{\operatorname{Ric}\}$. 
Since $U$ is parallel, we have that
$\mathcal{H}_f(X)=\nabla_X \nabla f=\lambda X$ if $X\in \mathcal{V}$ and that
$\mathcal{H}_f(X)=\nabla_X \nabla f=0$ if 
$X\in\ker\{\mathcal{H}_f\}=\operatorname{Image}\{\operatorname{Ric}\}$.
Consequently, the distribution $\mathcal{V}$ is parallel. 
Since the metric is not degenerate on $\mathcal{V}$, 
this implies that the manifold locally decomposes as a 
product $B\times F$ so that $B$ is Ricci flat and hence flat. On the other hand $F$ is 
Einstein satisfying $\rho^F=\lambda g^F$. 
Therefore the soliton is rigid. This completes the
proof of Theorem~\ref{T4}.\hfill\qed

\subsection{The proof of Theorem~\ref{T5}}\label{S3.3}
If $\dim(M)= 3$ the result follows from the discussion above since $\dim(\ker\{\operatorname{Ric}\})=2$. 
Assume $\dim(M)= 4$ henceforth. Using the previous discussion, we
need only examine the case  $\dim(\ker\{\operatorname{Ric}\})=3$. 
We are going to use Theorem~\ref{T4} to show that
 $\operatorname{Image}\{\operatorname{Ric}\}$ is a non-null parallel distribution. 
 We consider the adapted basis $\{U,\nabla f, E, F\}$ 
 where $\{U,\nabla f,E\}$ is a basis of 
 $\ker\{\operatorname{Ric}\}$ and 
 $F\cdot\mathbb{R}=\operatorname{Image}\{\operatorname{Ric}\}$.
 We show that the Weyl tensor is harmonic and $(M,g,f)$ is rigid by examining the
 components of the curvature tensor which have $\nabla f$ as an argument:
\medbreak\quad
$R(E,\nabla f,E,\nabla f)=(\nabla_{E} \rho)(\nabla f,E)-(\nabla_{\nabla f} \rho)(E,E)=0$,
\smallbreak\quad
$R(F,\nabla f,F,\nabla f)=(\nabla_{F} \rho)(\nabla f,F)-(\nabla_{\nabla f} \rho)(F,F)=0$,
\smallbreak\quad
$R(F,\nabla f,E,\nabla f)=\rho(F,E)\|\nabla f\|^2=0$,\qquad
\smallbreak\quad
$R(F,E, F,\nabla f)=\rho(\nabla f,E)=0$,\quad$R(E,F, E,\nabla f)= \rho(\nabla f,F)=0$.
\hfill\qed

\section{Steady locally homogeneous Lorentzian gradient Ricci solitons\\
The proof of Theorems~\ref{T7}--\ref{T8}}\label{S4}
Again, we shall use Lemma~\ref{L2} throughout the section without further citation.
Let $(M,g,f)$ be a steady locally homogeneous Lorentzian gradient Ricci soliton. 
Then $\|\operatorname{Hess}_f\|^2=0$ and
$\|\nabla f\|^2=\mu$ is constant.
In what follows we will consider the possibilities $\mu<0$ and $\mu=0$ separately.

\subsection{The proof of Theorem~\ref{T7}}
Assume that $\mu<0$. 
As $\mathcal{H}_f(\nabla f)=0$, we may restrict $\mathcal{H}_f$
to $\nabla f^\perp$. As $\nabla f^\perp$ inherits a positive definite metric and since 
$\|\operatorname{Hess}_f\|^2=0$, $\mathcal{H}_f=0$.
This shows that $\nabla f$ is a parallel vector field, and 
thus $(M,g)$ is locally a product 
$(\mathbb{R}\times N, -dt^2+g_N)$, where $(N,g_N)$ is a 
locally homogeneous Riemannian manifold  (see, for example, \cite{GK}). 
Additionally, $(N,g_N)$ is a steady gradient Ricci soliton,
and therefore Ricci flat. Following \cite{Spiro}, locally homogeneous Ricci flat Riemannian manifolds are locally isometric to Euclidean space. This
completes the proof of Theorem~\ref{T7}.

\subsection{The proof of Theorem~\ref{T8}~(1)}
Assume that $\|\nabla f\|^2=0$ so $\nabla f$ is a null vector. 
Choose an orthonormal
basis $\{E_1,...,E_{n+2}\}$
for the tangent space at a point so $E_1$ is timelike, so 
$\{E_2,...,E_{n+2}\}$ are spacelike,
and so $\nabla f=c(E_1+E_2)$ for some $c\ne0$. We further normalize the basis so
$\mathcal{H}_fE_1\in\operatorname{Span}\{E_1,E_2,E_3\}$. 
Let $\mathcal{H}_fE_i=\mathcal{H}_i^jE_j$.
Since $E_1+E_2\in\ker\{\mathcal{H}_f\}$, 
$\mathcal{H}_1^i+\mathcal{H}_2^i=0$ for all $i$.
Furthermore,  
$\mathcal{H}_1^i=\mathcal{H}_2^i=0$ for $i\ge4$
since $\mathcal{H}_fE_1\in\operatorname{Span}\{E_1,E_2,E_3\}$. 
Finally, since $\mathcal{H}_f$ is self-adjoint,
$\mathcal{H}_1^i=-\mathcal{H}_i^1$ for $2\le i$ and 
$\mathcal{H}_i^j=\mathcal{H}_j^i$ for 
$2\le i,j$. We summarize these relations:
\begin{equation}\label{E8}
\begin{array}{ll}
\mathcal{H}_1^i=-\mathcal{H}_i^1\text{ for }i\ge2,&
\mathcal{H}_i^j=\mathcal{H}_j^i\text{ for }2\le i,j,\\
\mathcal{H}^i_1=\mathcal{H}^i_2=0\text{ for }i\ge 4,& \mathcal{H}^i_1+\mathcal{H}^i_2=0\text{ for all }i
\vphantom{\vrule height 11pt}.
\end{array}\end{equation}
Since $\mathcal{H}_f=\mathcal{H}_i^jE^i\otimes E_j$
and $\|\operatorname{Hess}_f\|^2=\lambda((n+2)\lambda-\tau)=0$, we have
\begin{equation}\label{E9}
0=\| \operatorname{Hess}_f\|^2=\|\mathcal{H}_f\|^2=
(\mathcal{H}_1^1)^2-2\sum_{i\ge2}(\mathcal{H}_i^1)^2
+\sum_{2\le j,k}(\mathcal{H}_j^k)^2\,.
\end{equation}
The relations of Equation~(\ref{E8}) then permit us to rewrite 
Equation~(\ref{E9}) in the form:
$$0=\sum_{3\le j,k}(\mathcal{H}_j^k)^2\,.$$
This implies $\mathcal{H}_j^k=0$ for $3\le j,k$ and thus by Equation~(\ref{E8}),
$\mathcal{H}_fE_i=0$ for $i\ge4$. 
Thus the relevant portion of the matrix $\mathcal{H}$ becomes:
$$\mathcal{H}=\left(\begin{array}{rrr}\mathcal{H}_1^1&
\mathcal{H}^1_2&\mathcal{H}^1_3\\
\mathcal{H}^2_1&\mathcal{H}^2_2&\mathcal{H}^2_3\\
\mathcal{H}^3_1&\mathcal{H}^3_2&\mathcal{H}^3_3
\end{array}\right)
= \left(\begin{array}{rrr}\mathcal{H}_1^1&-\mathcal{H}_1^1&\mathcal{H}^1_3\\
\mathcal{H}_1^1&-\mathcal{H}_1^1&\mathcal{H}^1_3\\
-\mathcal{H}^1_3&\mathcal{H}^1_3&0
\end{array}\right)\,.
$$
We compute
$$\mathcal{H}^2=(\mathcal{H}_1^3)^2
\left(\begin{array}{rrr}-1&1&0\\
-1&1&0\\0&0&0
\end{array}\right)\text{ and }\mathcal{H}^3=0\,.$$
This shows that $\mathcal{H}$ is either $2$ or $3$-step nilpotent which
proves Assertion~(1).

\subsection{The proof of Theorem~\ref{T8}~(2)}
 Let $\mathcal{H}_f$ be $2$-step nilpotent.
The analysis above shows $\nabla f\in\operatorname{Image}\{\mathcal{H}_f\}$.
Since $\mathcal{H}_f$ has rank $1$, 
$\operatorname{Image}\{\mathcal{H}_f\}=\nabla f\cdot\mathbb{R}$. 
We use the Fredholm alternative and the fact that $\mathcal{H}_f$
is self-adjoint to establish Assertion~(2a) using the following equivalencies:
$$\begin{array}{llll}
&\mathcal{H}_fZ=0&\Leftrightarrow&g(\mathcal{H}_fZ,Y)=0\ \forall\ Y\\
\Leftrightarrow&g(Z,\mathcal{H}_fY)=0\ \forall\ Y&
\Leftrightarrow&Z\perp\operatorname{Range}\{\mathcal{H}_f\}\\
\Leftrightarrow&Z\perp\nabla f\,.
\end{array}$$
Choose a vector field $U$ so $g(U,\nabla f)=1$.
Since $\operatorname{Range}\{\mathcal{H}_f\}=\nabla f$ and since 
$g(U,\nabla f)=1$,
the fact that $\nabla f$ is recurrent follows from the
equation:
\begin{equation}\label{E10}
\nabla_X (\nabla f)=\mathcal{H}_f(X)=\theta(X)\cdot\nabla f
\text{ where }\theta(X)=g(U,\mathcal{H}_f(X))\,.
\end{equation}
Let $X$ and $Y$ be smooth vector fields in $\nabla f^\perp$. We show
that $[X,Y]$ belongs to $\nabla f^\perp$ and thus $\nabla f^\perp$ is an
integrable distribution by computing:
\begin{eqnarray*}
&&g([X,Y],\nabla f)=g(\nabla_XY-\nabla_YX,\nabla f)\\
&=&Xg(Y,\nabla f)-g(Y,\nabla_X\nabla f)-Yg(X,\nabla f)+g(X,\nabla_Y  \nabla f)\\
&=&X\{0\}- \operatorname{Hess}_f(Y,X)-Y\{0\}+ \operatorname{Hess}_f(X,Y)=0\,.
\end{eqnarray*}
Let $\gamma(t)$ be a geodesic with $\dot\gamma(0)\perp\nabla f$. 
We compute
$$\partial_tg(\dot\gamma,\nabla f)= g(\ddot\gamma,\nabla f)+g(\dot\gamma,
\nabla_{\partial_t}\nabla f)=\theta(\partial_t)g(\dot\gamma,\nabla f)\,.$$
Since $g(\dot\gamma,\nabla f)(0)=0$, the fundamental theorem of ODE's implies
$g(\dot\gamma,\nabla f)$ vanishes identically and thus 
$\dot\gamma\in\nabla f^\perp$. Since $g(\dot\gamma,\nabla f)=\partial_tf$,
the geodesic lies entirely in the level set of $f$. Assertion~(2b) follows.

We proceed by induction on the dimension to establish Assertion~(2c). 
Fix a point $P\in M$.
Let $\mathcal{V}:=\operatorname{Span}\{U,\nabla f\}$. The metric on $\mathcal{V}$
is non-degenerate and contains a null vector; consequently $\mathcal{V}$ has
Lorentzian signature.
We can choose complementary
Killing vector fields $\{F_1,\dots,F_n\}$ so $\{U,\nabla f,F_1,\dots, F_n\}$ is a
local frame field near $P$ and so that 
\begin{equation}\label{E11}
g(U,F_i)|_P=g( \nabla f,F_i)|_P=0\,.
\end{equation}
Consequently $\operatorname{Span}\{F_1,\dots,F_n\}$ is spacelike near $P$. 
Let $\xi_i:=\grad\{F_i(f)\}$; these are parallel vector fields by Lemma~\ref{L2}.
Let $\mathcal{W}:=\operatorname{Span}\{\xi_1,\dots,\xi_n\}$. Since the $\xi_i$ are
parallel, $r(x):=\operatorname{Rank}\{\mathcal{W}(x)\}$ is locally constant. Suppose $r>0$.
By reordering the collection $\{F_1,\dots,F_n\}$ if necessary, we may assume that 
$\{\xi_1,\dots,\xi_r\}$ is a local frame field for $\mathcal{W}$. 
Let $\epsilon_{ij}:=g(\xi_i,\xi_j)$ describe the induced metric on $\mathcal{W}$. 
Again we use the fact that the $\xi_i$ are parallel; this implies that the
$\epsilon_{ij}$ are constant. We can diagonalize $\epsilon$ or equivalently renormalize
the choice of the Killing vector fields $F_i$ to assume that $\epsilon$ is in fact diagonal. 
If $\operatorname{det}(\epsilon)=0$, then $\xi_i$ is a parallel null vector field for some $i$ and  Assertion~(2-c-i) holds. Thus we may assume that
the inner-product restricted to $\mathcal{W}$ is non-degenerate.
We may use Theorem~\ref{T1}
to decompose, at least locally, $M=N^{2+n-r}\times{\mathbb{R}^r_\nu}$. If the metric on $N$
is Riemannian, we may apply Theorem~\ref{T1} to see that the
soliton is trivial. Thus $N$ is Lorentzian. If
$\dim(N)=2$, then Theorem~\ref{T6} shows $N$ is flat and
$\mathcal{H}_f=0$ which is false. This shows $\dim(N)\ge3$ and we
may use our induction hypothesis on $N$. Thus we may assume without
loss of generality that $r=0$ so $\mathcal{W}=\{0\}$ and assume
henceforth that:
\begin{equation}\label{E12}
\grad\{F_i(f)\}=0\text{ for all }i\,.\end{equation}
By Equation~(\ref{E12}), $\kappa_i:=F_i(f)$
is constant for all $i$. By Equation~(\ref{E11}), 
$$\kappa_i=F_i(f)|_P=g(F_i,\nabla f)|_P=0\,.$$
Consequently $g(F_i,\nabla f)$ vanishes identically and we have
\begin{equation}\label{E13}
F_i\in\ker\{\mathcal{H}_f\}=\ker\{\operatorname{Ric}\}=\nabla f^\perp\,.
\end{equation}
We may use Equation~(\ref{E10}) and Equation~(\ref{E13}) to see
\begin{equation}\label{E14}
\begin{array}{l}
\nabla_{\nabla f}\nabla f=\mathcal{H}_f(\nabla f)=0,\quad
\nabla_{F_i}\nabla f=\mathcal{H}_f( F_i)=0\text{ for all i},\\
\nabla_U\nabla f=\mathcal{H}_f (U)=
\Xi\nabla f\text{ where }\Xi:=g(\mathcal{H}_f (U),U)=-\rho(U,U)\,.
\vphantom{\vrule height 11pt}\end{array}\end{equation}
We use Equation~(\ref{E14}) to see:
\begin{equation}\label{E15}
\nabla_Y{ \nabla f}=0\text{ if }Y\perp\nabla f\,.
\end{equation}
Thus the only covariant derivative at issue is $\nabla_U{ \nabla f}$.
We shall let $\Psi:=\psi\cdot\nabla f$. This is a null vector field. By
Equation~(\ref{E15}), $\Psi$ will be parallel if and only if $\psi$ satisfies
the equations:
\begin{equation}\label{E16}
Y(\psi)=0\text{ if }Y\perp\nabla f\text{ and }U(\psi)+\psi\Xi=0\,.
\end{equation}
   
Since $F_i$
is a Killing vector field, $\nabla_{F_i}\rho=0$. Since
$F_i\in\ker\{\operatorname{Ric}\}$, $\rho(F_i,\cdot)$ vanishes identically.
Consequently, Lemma~\ref{L2} yields
\begin{equation}\label{E17}
\begin{array}{l}
R(F_i,U,F_j,\nabla f)=(\nabla_{F_i}\rho)(U,F_j)-(\nabla_U\rho)(F_i,F_j)\\
\quad=- U\rho(F_i,F_i)+\rho(\nabla_UF_i,F_j)+\rho(\nabla_UF_j,F_i)=0\,.
\vphantom{\vrule height 11pt}\end{array}\end{equation}
Let $g_{ij}=g(F_i,F_j)$.
Since $U\in\ker\{\operatorname{Ric}\}$, since $\{U,{ \nabla f}\}$ span a hyperbolic pair,
Equation~(\ref{E17}) implies:
\begin{eqnarray*}
0&=&\rho(U,\nabla f)|_{P}=R(U,\nabla f,\nabla f,U)|_{P}
+\displaystyle\sum_{i,j=1}^n g^{ij}R(U,F_i,\nabla f,F_j)|_{P}\\
&=&R(U,\nabla f,\nabla f,U)|_{P}\,.
\end{eqnarray*}
Since $P$ was arbitrary and the only condition on $U$ was that $g(U,{ \nabla f})=1$,
this holds for arbitrary $P$ and we have
\begin{equation}\label{E18}
0=R(U,\nabla f,\nabla f,U)\text{ if }g(U,\nabla f)=1\,.
\end{equation}
Also, in general, if $X$ is a Killing vector field, then for arbitrary
vector fields, we have (see, for example, \cite{Kostant, Nomizu1}) that:
$$R(X,Y)Z=-\nabla_Y\nabla_ZX+\nabla_{\nabla_Y Z}X\,.$$
Let $\Xi$ be as defined in Equation(\ref{E14}).
We use Equation~(\ref{E13}) to see:
$$g(\nabla_UF_i,\nabla f)=U\,g(F_i,\nabla f)-g(F_i,\nabla_U{ \nabla f})
=-g(F_i,\Xi \nabla f)=0\,.$$
Since the $F_i$ are Killing vector fields, since $g(F_i,\nabla f)=0$,
and since $\nabla f$ is
recurrent, 
\begin{eqnarray*}
&&R(F_i,U,U,\nabla f)=-g(\nabla_U\nabla_UF_i,\nabla f)
+g(\nabla_{\nabla_UU}F_i,\nabla f)\\
&=&-U\,g(\nabla_UF_i,\nabla f)+g(\nabla_UF_i,\nabla_U{\nabla f})
\\&&\quad+(\nabla_UU)g(F_i,\nabla f)
-g(F_i,\nabla_{\nabla_UU}\{\nabla f\})\\
&=& -U\{U\,g(F_i,\nabla f)-g(F_i,\nabla_U{\nabla f})\}+g(\nabla_UF_i,\Xi \nabla f)\\
&=&U\, g(F_i,\Xi \nabla f)+\Xi g(\nabla_UF_i,\nabla f)=0\,.
\end{eqnarray*}

By Lemma~\ref{L2}, if $\{X\,,Y\,,Z\}$ are vector fields on
a gradient Ricci soliton, then
$$R(X,Y,Z,\nabla f)=(\nabla_X \rho)(Y,Z)-(\nabla_Y \rho)(X,Z)\,.$$
Consequently, we have that
\begin{eqnarray*}
&&0=R(U,\nabla f,U,\nabla f)=(\nabla_U \rho)(\nabla f,U)-(\nabla_{\nabla f} \rho)(U,U)\\
&&0=R(F_i,U,U,\nabla f)=(\nabla_{F_i} \rho)(U,U)-(\nabla_U \rho)(F_i,U).
\end{eqnarray*}
By Equation~(\ref{E14}), $\Xi=-\rho(U,U)$. Thus we may compute:
\medbreak\quad
$-\nabla f(\Xi)=\nabla f \rho(U,U)=(\nabla_{\nabla f} \rho)(U,U)
+2\,\rho(\nabla_{\nabla f}U,U)$
\smallbreak\qquad
$= (\nabla_U \rho)(\nabla f,U)-2\, g(\nabla_{\nabla f}U,\Xi \nabla f)$
\smallbreak\qquad
$= U\,\rho(\nabla f,U)-\rho(\nabla_U\nabla f,U)-\rho(\nabla f,\nabla_UU)$
\smallbreak\qquad\quad
$-2\Xi(\nabla fg(U,\nabla f)-g(U,\nabla_{\nabla f}\nabla f))=0$, and
\medbreak\quad
$-F_i(\Xi)=F_i \rho(U,U)
= (\nabla_{F_i} \rho)(U,U)+2\,\rho(\nabla_{F_i} U,U)$
\smallbreak\qquad
$= (\nabla_U \rho)(F_i,U)-2\, g(\nabla_{F_i} U, \Xi \nabla f)$
\smallbreak\qquad
$= U\,\rho(F_i,U)-\rho(\nabla_U F_i,U)-\rho(F_i,\nabla_UU)$
\smallbreak\qquad\quad
$-2\Xi(F_ig(U, \nabla f)-g(U,\nabla_{F_i} \nabla f))$
\smallbreak\qquad
$= g(\nabla_U F_i,\Xi \nabla f)=\Xi U g(F_i,\nabla f)-\Xi g(F_i,\Xi \nabla f)=0$.
\medbreak\noindent
This shows that $X(\Xi)=0$ if $X\in\nabla f^\perp$. Since the distribution
$\nabla f^\perp$ is integrable, the Frobenius theorem means we can introduce
local coordinates $(u,x^2,...,x^{n+2})$ so that $U=\partial_u$ and
$\nabla f^\perp=\operatorname{Span}\{\partial_{x_2},...,\partial_{x_{n+2}}\}$.
Thus Equation~(\ref{E16}) becomes an ordinary differential equation which
can be solved. This completes the proof of Theorem~\ref{T8}.\hfill\qed

\begin{example}\label{Exm1}
\rm We follow the discussion in \cite{BBGG}. A Cahen-Wallach space has
the metric  given locally by Equation~(\ref{E4}):
$$
g=2\, dt dy+\left(\sum_{i=1}^{n} \kappa_i\, x_i^2\right) dy^2+\sum_{i=1}^{n} dx_i^2
\text{ for }0\ne\kappa_i\in\mathbb{R}\,.
$$
 The Levi-Civita connection
is determined by the non-zero Christoffel symbols:
$$\nabla_{\partial_y}\partial_y=-\sum_i\kappa_i x_i\partial_{x_i}\text{ and }
\nabla_{\partial_y}\partial_{x_i}=\nabla_{\partial_{x_i}}\partial_y=\kappa_i x_i\partial_v\,.$$
Thus the only non-zero entries in  the curvature tensor are given by:
$$R(\partial_y, \partial_{x_i},\partial_y, \partial_{x_i})=-\kappa_i$$
and thus (possibly) non-zero entries in the Ricci tensor are 
$$\rho(\partial_y,\partial_y)=-\kappa\text{ where }\kappa:=\kappa_1+...+\kappa_n\,.$$
Assuming that $\kappa\ne0$, we then have
$\operatorname{Ric}{(\partial_y)}=-\kappa\partial_t$ and 
$\operatorname{Ric}{(\partial_t})=0$.
Thus the Ricci tensor is two step nilpotent. The $f$ defines
a gradient Ricci soliton if and only if $f(t,y,x_1,...,x_n)=f(y)$ where 
$f(y)=a_0+a_1y+\frac14\kappa y^2$; $\lambda=0$ in this instance.
Note that $df= (a_1+\frac12\kappa y) dy$ and hence $\nabla f= (a_1+\frac12\kappa y)\partial_t$
is a null parallel vector field.
\end{example}\section{Symmetric gradient Ricci solitons\\
The proof of Theorem~\ref{T11}}\label{S5}

Let $(M,g)$ be a locally symmetric Lorentzian 
manifold. If $(M,g,f)$ is a  non-steady gradient Ricci soliton, then by
Theorem~\ref{T3}, $M$ splits, at least locally,
 as a product  $M=N_0\times N_1\times\mathbb{R}^k$, where $(N_0,g_0)$ is indecomposable but reducible and $(N_1,g_1)$ is Einstein. If $N_0$ does not appear in the decomposition, then the soliton is rigid. Otherwise, $(N_0,g_0)$ is an indecomposable but not irreducible Lorentzian 
symmetric space, hence a Cahen-Wallach symmetric space \cite{CW70} (see also \cite{BI}). Theorem~\ref{T10}
rules out this latter possibility since if $(N,g_N,f_N)$ is a Cahen-Wallach
gradient Ricci soliton, then it is steady.

Next suppose that $(M,g,f)$ is a locally symmetric Lorentzian  steady gradient Ricci soliton.
We can use the 
de Rham-Wu decomposition of the manifold to split $(M,g)$ locally as a product
${M}= N\times M_1\times\cdots\times M_l\times \mathbb{R}^k_\nu$, where $(N,g_N)$ is a Cahen-Wallach
symmetric space, where the $M_i$ are irreducible symmetric spaces, and where $\mathbb{R}^k_\nu$ is either Euclidean or Minkowskian space.
Since irreducible symmetric spaces are Einstein, the induced soliton is either trivial 
or the scalar  curvature vanishes, which implies that $M_i$ is Ricci flat.
If $M_i$ is Riemannian, then it is flat since Ricci flat 
locally symmetric spaces are flat in the 
Riemannian setting \cite{besse,helgason}. 
Moreover, if $M_i$ is Lorentzian, then it is flat since irreducible Lorentzian 
locally symmetric spaces are of constant sectional curvature \cite{CLPTV1}.
Hence, if the gradient Ricci soliton is steady, 
then the decomposition above reduces to 
${M}=N\times\mathbb{R}^k$, where $(N,g_N)$ is a Cahen-Wallach symmetric space. 
Theorem~\ref{T11} now follows.\hfill\qed

\section{Three-dimensional locally homogeneous gradient Ricci solitons}\label{S6}

\subsection{The proof of Theorem~\ref{T12}}
Let $(M,g)$ be a 3-dimensional Lorentzian  strict Walker metric. 
There exist local coordinates so the metric is given by  Equation~\eqref{E5x}: 
$$g= 2 dt dy + dx^2+\phi(x,y)dy^2\,.
$$
Let $f(t,x,y)$ be a smooth real valued function. To simplify the notation,
set $f_t=\frac{\partial f}{\partial t}$, 
  $f_{tx}=\frac{\partial^2 f}{\partial t \partial x}$, and so forth.
One computes easily that the soliton equation 
$\operatorname{Hess}_f+\rho=\lambda g$ is equivalent to the following relations:
\begin{equation}\begin{array}{ll}\label{E19}
0=f_{tt}=f_{tx},&
0=f_{xx}-\lambda=f_{ty}-\lambda,\\
0=2f_{xy}-\phi_{x}f_t,&
0=2\lambda\,\phi+\phi_{xx}-2 f_{yy}-\phi_x f_x+\phi_y f_t.
\end{array}\end{equation}
We use the first identities in Equation~(\ref{E19}) to see:
\[
f(t,x,y)=t(\lambda\, y +\kappa)
+\textstyle\frac{1}{2}\,\lambda\, x^2+  \alpha(y)\,x +\gamma(y)\text{ for }\kappa\in\mathbb{R}\,.
\]
Hence, the equations of Equation~(\ref{E19}) simplify to become:
\begin{eqnarray}
&&0=2\,\alpha'(y)-\left(\lambda\, y+\kappa\right)\,\phi_x,\label{h33}\\
&&0=2\,\lambda\,\phi-2\,\gamma''(y)-2\,x\, \alpha''(y)+(\lambda\, y+\kappa)\,\phi_y
- (\lambda\, x + \alpha(y))\,\phi_x+\phi_{xx}.\label{h44}
\end{eqnarray}
We differentiate Equation~\eqref{h33} with respect to $x$  to conclude:
\begin{equation}\label{h5}
0=(\lambda\, y + \kappa)\,\phi_{xx}.
\end{equation}
Since the Ricci operator is given by:
\[
\operatorname{Ric}=\left(\begin{array}{ccc}
0&0&-\frac{1}{2}\,\phi_{xx}\\
0&0&0\\
0&0&0
\end{array}\right),
\]
the metric is flat if and only if $\phi_{xx}=0$. 
Since we assume that the Walker metric is not-flat, we may use 
Equation~(\ref{h5}) to see that $\lambda=\kappa=0$ and conclude that the
gradient Ricci soliton is steady. Consequently Equations \eqref{h33} and \eqref{h5}
imply that $f(t,x,y)= \alpha\,x+\gamma(y)$ so Equation~\eqref{h44} becomes:
\begin{equation}\label{h6}
2\,\gamma''(y)+\alpha\,\phi_x-\phi_{xx}=0.
\end{equation}
We take the derivative with respect to $x$ to see
$\alpha\,\phi_{xx}=\phi_{xxx}$. We examine the two cases seriatim.
\subsection*{Case I: Suppose that $\alpha\ne0$}
We then have:
$$
\phi(x,y)=\frac{1}{\alpha^2}\,a(y)\,e^{\alpha x}+x\, b(y)+c(y)
$$
for some arbitrary functions $a(y)\neq 0$, $b(y)$ and $c(y)$.
Moreover the potential function of the soliton is given by 
$f(t,x,y)$ $=$ $\alpha\, x +\gamma(y)$, where $\gamma''(y)=-\frac{1}{2}\, \alpha\, b(y)$.
In this case $\nabla f=\gamma'(y)\,\partial_t+\alpha\,\partial_x$ is spacelike.
This gives rise to the first possibility in Theorem~\ref{T12}.
\subsection*{Case II: Suppose that $\alpha=0$} We then have:
$$
\phi(x,y)=x^2\, a(y)+x\, b(y)+ c(y)
$$
for some arbitrary functions $a(y)\neq 0$, $b(y)$ and $c(y)$. 
Moreover the potential function of the soliton is given by 
$f(t,x,y)=\gamma(y)$, where $\gamma''(y)=\frac14\, a(y)$.
In this case $\nabla f=\gamma'(y)\,\partial_t$ is a null and recurrent vector field.
This gives rise to the second possibility in Theorem~\ref{T12}.

\subsection{The proof of Theorem~\ref{T15}}
Let $(M,g,f)$ be a locally homogeneous Lorentzian gradient Ricci soliton of dimension 3. 
\subsection*{Case I: Suppose that $(M,g,f)$ is non steady}
By Theorem~\ref{T5} the soliton is rigid. 

\subsection*{Case II: Suppose that $(M,g,f)$ is steady} Consequently by Lemma~\ref{L2}, 
the potential function is a  solution of the Eikonal equation  $||\nabla f||^2=\mu$. We distinguish
3 subcases:
\subsection*{Case II-a: $(M,g)$ is steady and $\mu<0$}
We apply Theorem~\ref{T7} to see that $(M,g)$ splits locally as a product
and hence the soliton is rigid.

\subsection*{Case II-b: $(M,g)$ is steady and $\mu=0$}
We use Theorem~\ref{T8} to see that the Ricci operator is either 2 or 3 step
nilpotent. It follows from work of \cite{CK} that there do not exist locally homogeneous $3$-dimensional
manifolds with $3$-step nilpotent Ricci operator. Consequently, the Ricci operator is $2$-step nilpotent
and $(M,g)$ admits a locally defined parallel null vector field by Theorem~\ref{T8}. Consequently, $(M,g)$
is locally a strict Walker manifold. Consequently, the underlying geometry of $(M,g)$ is given by Theorem~\ref{T14};
the function $f$ is now determined by Theorem~\ref{T12}.

\subsection*{Case II-c: $(M,g)$ is steady and $\mu>0$}
Since the scalar curvature is constant, the Ricci operator satisfies $\operatorname{Ric}(\nabla f)=0$, 
which shows that either $f$ is constant, or otherwise the Ricci operator has a zero eigenvalue. 
We now consider the different possibilities for the kernel of $\operatorname{Ric}$.

Assume $\dim(\ker\{\operatorname{Ric}\})=1$. 
It follows from \cite{Calvaruso1}  that
$(M,g)$ is either a symmetric space or a Lie group. If $(M,g)$ is symmetric, 
then it is one of the following: a manifold of constant sectional curvature, a product 
$\mathbb{R}\times N$ where $(N,g_N)$ is of constant curvature, or a three-dimensional 
Cahen-Wallach symmetric space. Hence, in all the cases, any gradient Ricci soliton is trivial, 
rigid or the underlying manifold admits a null parallel vector field (and we have already examined
that case). Now we concentrate on Lie groups. 
Since the eigenspaces of the Ricci operator are left-invariant, since $\nabla f$ has constant norm $\mu>0$, 
and since $\dim(\ker\{\operatorname{Ric}\})=1$ we have that $\nabla f$ is a 
left-invariant vector field.
Left-invariant Ricci solitons on three-dimensional Lorentzian Lie groups were considered in \cite{israel},
showing that they exist in and only if the Ricci operator has exactly one-single eigenvalue, 
which must be zero since $\operatorname{Ric}(\nabla f)=0$.
This shows that the Ricci operator is three-step nilpotent, but that is not possible due to the analysis carried out in \cite{CK}.

Finally assume $\dim(\ker\{\operatorname{Ric}\})=2$. 
In this case the Ricci operator is either diagonalizable or two-step nilpotent. 
The later implies  that the manifold admits locally a null parallel vector field \cite{CGVV}, 
and again this case has been treated. If the Ricci operator is diagonalizable, then 
$\|\operatorname{Ric}\|^2=\pm \tau^2=\|\operatorname{Hess}_f\|^2$ and Lemma~\ref{L2}~(3) 
shows that $\tau=0$, from where it follows that $(M,g)$ is flat and the soliton is trivial. This completes
the proof of Theorem~\ref{T15}.
\hfill\qed

\end{document}